\documentclass[11pt]{article}
\usepackage{amsfonts}
\usepackage{amsthm}
\usepackage{amsmath}
\usepackage{amssymb}
\usepackage{amscd}
\usepackage{verbatim}
\usepackage{multicol}

\theoremstyle{definition}
\newtheorem{thm}{Theorem}[section]
\newtheorem{prop}[thm]{Proposition}

\numberwithin{equation}{section}

\def\cqd{\rightline{$\qedsymbol$}}
\def\obs{\noindent{\normalsize{\bf Remark: }}}
\def\bb#1{\text{$\mathbb{#1}$}}
\def\cal#1{\text{$\mathcal{#1}$}}
\def\so{\Rightarrow}

\def\lie#1{\mathfrak{#1}}
\def\tlie#1{\widetilde{\mathfrak{#1}}}
\def\hlie#1{\widehat{\mathfrak{#1}}}
\def\kz{Knizhnik-Zamolodchikov }

\begin{document}

\author{Adriano Adrega de Moura\footnote{Ph.D. studies supported by FAPESP (99/11600-0) and visit to MIT supported by CAPES (0365/01-7), Brazil}\\ {\small adrianoam@ime.unicamp.br}}
\title{Elliptic Dynamical R-Matrices from the\\ Monodromy of the q-\kz Equations\\ for the Standard Representation of $U_q(\tlie{sl}_{n+1})$}
\date{December/2001}
\maketitle


\section*{Introduction}

\nocite{efk01}
\nocite{eo01}
\nocite{ev97}
\nocite{ev99}
\nocite{fel94}
\nocite{felv99}
\nocite{fr01}
\nocite{gara01}
\nocite{jim02}
\nocite{tv97}
\nocite{ww01}

In 1984 Knizhnik and Zamolodchikov derived differential equations for conformal blocks for the Wess-Zumino-Witten conformal field theory. These remarkable equations are now known as the Knizhnik-Zamolodchikov (KZ) equations. 

A few years later Kohno and Drinfeld considered the monodromy of the KZ equations, and showed that it is given by R-matrices for the corresponding quantum group. This theorem is called the Drinfeld-Kohno theorem. 

In 1992, in the pioneering work \cite{fr01}, Frenkel and Reshetikhin generalized the notion of conformal blocks to the q-deformed case, and proposed a q-deformed version of the Knizhnik-Zamolodchikov equations -- the so called qKZ equations, which are no longer differential but rather q-difference equations. Furthermore, they calculated the connection (= q-monodromy) of these equations in a particular case, and showed that they coincide with certain elliptic solutions of the star-triangle relations known in statistical mechanics. On the basis of these calculations, they conjectured that the Drinfeld-Kohno theorem for KZ equations should have a q-analog, saying that connection matrices of qKZ equations should be R-matrices of an appropriate elliptic quantum group. This statement had to remain a conjecture, since at that time there was no sensible definition of elliptic quantum groups. 

A successful theory of elliptic quantum groups began with the work of Felder \cite{fel94}. In this work it was realized for the first time that the main equation in this theory is not the usual quantum Yang-Baxter equation, but rather the dynamical quantum Yang-Baxter equation (DQYBE), which appeared earlier in physics literature. In particular, Felder defined elliptic quantum groups of type A, using the corresponding solutions of the DQYBE. 

After this development, it became possible to formulate the Frenkel-Reshetikhin conjecture precisely (at least in type A) : the connection matrices for qKZ equations are R-matrices for Felder's elliptic quantum groups. This was proved by Tarasov and Varchenko first for $\lie{sl}_2$ \cite{tv97} and then for $\lie{sl}_{n+1}$ (to be published), using a rather sophisticated theory of integral formulas. 

A different approach to this conjecture became possible after works of Etingof and Varchenko \cite{ev99}, where the so-called exchange construction was developed. Using this approach, it is possible to derive the general formulation of the Frenkel-Reshetikhin conjecture from its version in a very simple special case -- in the tensor product of two copies of the vector representation of $\lie{sl}_{n+1}$. This special case can be handled by means of the classical theory of the q-hypergeometric function $_2\phi_1$, exactly as it was done in \cite{fr01}. 

The development of this special case in full detail is the main goal of this paper. It seems that it has not been done previously, although the ideas behind it have certainly been known for a long time. We plan to consider the general case and other applications in a separate publication. 

The structure of this paper is as follows. 

In section 1, we fix the basic notation, recall some definitions and known facts, and state the main results, which are the formula for the
exchange matrix and the sequence of gauge transformations \cite{ev97} that transforms it into Felder's solution.We note that one of the gauge transformations needed is multiplication by a scalar function which is a product of elliptic gamma functions of 
Felder and Varchenko \cite{felv99}; the meaning of its appearance in this context is not clear to us.
 
The other sections contain the proofs of these results. Namely, in Section 2 we calculate the action of the universal R-matrix of $U_q(\lie{sl}_{n+1})$ on the product of vector representations and solve the q-\kz equations, while in Section 3 we obtain the elliptic dynamical R-matrix through the exchange construction. Some of the intermediate proofs, which are straightforward, but lengthy, are in the appendix. We also included in the appendix the standard facts and formulas which are used througout the paper. We note that while performing the calculations we found a number of computational errors in some previous publications. They are listed in the last section of the appendix.

{\bf Acknowledgements.} The problem solved in this paper was suggested to me by professor P. Etingof during my visit to MIT as a Ph.D student. It is motivated by his work with A. Varchenko \cite{ev97,ev99}. I would like to thank professor Etingof, not only for suggesting the problem, but also for all helpful explanations.


\section{The Main Theorems}

\subsection{Recalling Definitions and Fixing Notation}

Consider $\lie g = \lie{sl}_{n+1}(\bb C)$ acting naturally on $V = \bb C^{n+1}$. Let $\{v_0, \dots , v_n\}$ be the canonical basis of $V$ and $E_{i,j}$ be the elementary matrix  having $0$ in all of its entries except in the $ij$-th which is $1$ and set $h_i = E_{i-1, i-1}-E_{i,i}, e_i = E_{i-1,i}, f_i = E_{i,i-1}$. Denote by $<,>$ the invariant bilinear form normalized that $<\alpha_i,\alpha_i> = 2$ for simple roots ($i = 1, \dots, n$). The dominant weights are $\{\omega_i = h_i^*\}$, the Coxeter number is $h\check{}=n+1$ and $\rho = \sum\omega_i$. The highest weight vector of V is $v_0$ and its weight is $\omega_1$, $v_n$ has weight $-\omega_n$ and the others $v_m$ have weight $\mu_m = \omega_{m+1} - \omega_m$. 

Sometimes it will be more convenient to identify $\lie h$ with a subspace of $\bb C^{n+1}$ setting $h_i = (0,0,\dots,1_{i-1},-1_i, 0,\dots,0) \in \bb C^{n+1}$. Then $<,>$ is just $< (\lambda_0, \dots, \lambda_n),(\mu_0, \dots, \mu_n)> = \sum \lambda_i\mu_i$ and we use it to identify $\lie h^*$ with $\lie h$. In appendix \ref{ap:useful} we collect some useful formulas like $<\mu_m,\mu_l>$ obtained easily from this identification.

Let $(a_{ij})$ be the extended Cartan matrix corresponding to $\widehat{\lie{sl}}_{n+1}$ where $c$ will denote the central element. Then the quantum group $U_q(\widetilde{\lie{sl}}_{n+1})$ is the algebra generated by $q^{td}, q^h, e_i, f_i,$ $i=0, \dots, n, t \in \bb C, h \in \hat{\lie h}$ with relations

$$q^{a+b} = q^a q^b, \quad a,b \in \tlie h = \hlie h + \bb Cd$$
$$q^h e_i q^{-h} = q^{\alpha_i(h)}e_i, \quad q^h f_i q^{-h} = q^{-\alpha_i(h)}f_i, \quad i=0,\dots, n, \, h \in \tlie h$$
$$e_if_j - f_je_i = \delta_{ij} \frac{q^{h_i}-q^{-h_i}}{q-q^{-1}}$$
$$\sum_{m=0}^{1-a_{ij}} \frac{(-1)^m}{[m]_q![1-a_{ij}-m]_q!} e_i^m e_j e_i^{1-a_{ij}-m} = 0$$
$$\sum_{m=0}^{1-a_{ij}} \frac{(-1)^m}{[m]_q![1-a_{ij}-m]_q!} f_i^m f_j f_i^{1-a_{ij}-m} = 0$$
$$\text{where} \quad [m]_q = \frac{q^m - q^{-m}}{q-q^{-1}}, \quad [m]_q! = [m]_q[m-1]_q\dots [1]_q$$

$U_q(\widetilde{\lie{sl}}_{n+1})$ is a quasi-triangular Hopf algebra whose structure maps are given by

$$\Delta(q^h) = q^h\otimes q^h, \quad \Delta(e_i) = e_i\otimes q^{h_i} + 1\otimes e_i, \quad \Delta(f_i) = f_i\otimes 1 + q^{-h_i}\otimes f_i$$
$$\varepsilon(q^h) = 1, \quad \varepsilon(e_i) = \varepsilon(f_i) = 0$$
$$S(q^h) = q^{-h}, \quad S(e_i) = -e_iq^{-h_i}, \quad S(f_i) = -q^{h_i}f_i$$

The universal R-matrix is of the form
$$\cal R = q^{c\otimes d + c\otimes c + \sum x_i\otimes x_i}(\sum a_j\otimes a^j)$$
where $\{x_i\}$ is an orthonormal basis for $\mathfrak h$, $\{a_j\}$ is a basis for $U_q(\mathfrak{\hat n}^+)$ and $\{a^j\}$ is the dual basis for $U_q(\mathfrak{\hat n}^-)$.


\subsection{Quantum Dynamical Yang-Baxter Equation with Spectral Parameter}

Let $\lie h$ be a finite dimensional abelian Lie algebra over \bb C, $V$ be a finite dimensional semisimple $\lie h$-module and $\gamma \in \bb C^*$. The quantum dynamical Yang-Baxter equation with spectral parameter and step $\gamma$ is 
\begin{align}\notag
R^{12}(u_1-u_2,\lambda-\gamma h^{(3)}) R^{13}(u_1-u_3,\lambda) R^{23}(u_2-u_3,\lambda-\gamma h^{(1)})\\ \label{eq:qdybs}
= R^{23}(u_2-u_3,\lambda) R^{13}(u_1-u_3,\lambda-\gamma h^{(2)}) R^{12}(u_1-u_2,\lambda)
\end{align}
with respect to a meromorphic function $R:\bb C\times\lie h^* \to \text{End}(V\otimes V)$. The notation $h^{(i)}$ is standard and means that, for example, $R^{12}(u_1-u_2,\lambda-\gamma h^{(3)})(v_1\otimes v_2\otimes v_3) = \big(R(u_1-u_2,\lambda-\gamma\mu)(v_1\otimes v_2)\big)\otimes v_3$ if $v_3$ has weight $\mu$. Zero-weight solutions of \eqref{eq:qdybs} are called quantum dynamical R-matrices.

Solutions of the form
$$R(u,\lambda) = \sum_{m,l}  \alpha_{m,l}(u,\lambda) E_{m,m}\otimes E_{l,l} + \sum_{m\neq l} \beta_{m,l}(u,\lambda) E_{l,m}\otimes E_{m,l}$$
are said to be of $\lie{gl}_{n+1}$ type. Felder \cite{fel94} found the interesting elliptic solutions
$$R_{\gamma,\tau}^{ell}(u,\lambda) = \sum_{m=0}^{n} E_{m,m}\otimes E_{m,m} + \sum_{m\neq l} \alpha(u,\lambda_{m,l}) E_{m,m}\otimes E_{l,l} + \beta(u,\lambda_{m,l}) E_{l,m}\otimes E_{m,l}$$
with
$$\alpha(u,\lambda) = \frac{\vartheta_1(\lambda+\gamma;\tau)}{\vartheta_1(\lambda;\tau)} \frac{\vartheta_1(u;\tau)}{\vartheta_1(u-\gamma;\tau)}\quad\quad \beta(u,\lambda) = \frac{\vartheta_1(\gamma;\tau)}{\vartheta_1(\lambda;\tau)} \frac{\vartheta_1(u-\lambda;\tau)}{\vartheta_1(u-\gamma;\tau)}$$
where $\lambda_{m,l} = \lambda_m - \lambda_l$, and $\vartheta_1$ is the standard first theta function \eqref{eq:ftheta}.

In \cite{ev97}, Etingof and Varchenko introduced gauge transformations of quantum dynamical R-matrices of $\lie{gl}_{n+1}$ type. We list three gauge transformations that will be important for us.
\begin{equation}\label{eq:g4}
R(u,\lambda) \mapsto c(u)R(u,\lambda)
\end{equation}
where $c(u)$ is a holomorphic scalar function;
\begin{equation}\label{eq:g5}
R(u,\lambda) \mapsto R(au,b\lambda + \mu)
\end{equation}
where $a,b \in \bb C^*$ and $\mu \in \lie h^*$;
\begin{gather}\notag
R(u,\lambda) \mapsto \\\label{eq:g2}
\sum_{m=0}^{n} \alpha_{m,m}(u,\lambda)E_{m,m}\otimes E_{m,m} + \sum_{m\neq l} \varphi_{m,l}(\lambda)\alpha_{m,l}(u,\lambda) E_{m,m}\otimes E_{l,l} + \beta_{m,l}(u,\lambda) E_{l,m}\otimes E_{m,l}
\end{gather}
where $\{\varphi_{m,l}(\lambda)\}$ is a $\gamma$-closed 2 form as defined in \cite{ev97}.

The only one that changes the step is the second. It sends $\gamma$ to $\gamma/b$.


\subsection{Fusion and Exchange Construction}\label{sec:fec}

One way of constructing solutions of \eqref{eq:qdybs} is through the so called exchange construction that we briefly recall now. Let $M_{\lambda}$ be the Verma module over $U_q(\lie g)$. For generic $\lambda$ and $q$ it is irreducible. Consider the induced $U_q(\tlie g)$-module 
$$M_{\lambda,k} = U_q(\tlie g)\otimes_{U_q(\tlie g^{\geq 0})}M_{\lambda}$$
where $U_q(\tlie g^{\geq 0})$ acts on $M_{\lambda}$ by $e_0 \equiv 0,\, q^c \equiv q^k$ for some scalar $k$ (the level of the representation) and $q^d \equiv q^{-\Delta}$ where
$$\Delta = \Delta_k(\lambda) = \frac{<\lambda,\lambda+2\rho>}{2\kappa} \quad\quad\quad \kappa = k+h\check{}$$
For generic $\lambda$ and $k$, $M_{\lambda,k}$ is irreducible as $U_q(\hlie g)$-module. Then we consider intertwining operators $\Phi(z) : M_{\lambda,k} \to M_{\mu,k}\hat{\otimes}V(z)$ where $V(z)$ is an evaluation representation of $U_q(\hlie g)$. We have the following characterization \cite{efk01,fr01} :

\begin{thm}\label{thm:intchar}
Let $M_{\lambda}, M_{\mu}$ be irreducible Verma modules over $U_q(\lie g)$ and $k\in\bb C^*$ such that $M_{\lambda,k}, M_{\mu,k}$ are irreducible. For every $U_q(\lie g)$-intertwiner $\varphi : M_{\lambda} \to M_{\mu}\otimes V$, there is a unique $U_q(\hlie g)$-intertwiner 
$$\Phi^{\varphi}_k(z) : M_{\lambda,k} \to M_{\mu,k}\hat{\otimes} V(z)$$
such that for any degree-zero vector $w$ (i.e., $w\in M_{\lambda_1}$), the degree-zero component of $\Phi^{\varphi}_k(z)w$ is $\varphi w$. The same holds with $V[z,z^{-1}]$.
\end{thm}

\begin{prop}\label{prop:tintchar}
The operator
$$\tilde{\Phi}^{\varphi}_k := z^{-\Delta}\Phi^{\varphi}_k(z) : M_{\lambda,k} \to M_{\mu,k}\hat{\otimes}z^{-\Delta}V[z,z^{-1}]$$
is a $U_q(\tlie g)$-intertwiner if and only if $\Delta = \Delta_k(\lambda) - \Delta_k(\mu)$.
\end{prop}  

Given an intertwiner $\varphi:M_{\lambda} \to M_{\mu}\otimes V$ where $V$ is a finite dimensional representation of $U_q(\lie g)$, the vector $<\varphi> \in V$ defined by
$$<\varphi> = <v_{\mu}^*,\varphi v_{\lambda}>$$
where $v_{\lambda}$ is the highest weight vector of $M_{\lambda}$ and $v_{\mu}^*$ is the lowest weight vector of the restricted dual $M_{\mu}^*$, is called the expectation value of $\varphi$. In the same lines of theorem \ref{thm:intchar} one shows that $<>$ is an isomorphism Hom$_{U_q(\lie g)}(M_{\lambda},M_{\mu}\otimes V) \to V[\lambda-\mu]$ (see \cite{eo01}). Then, given a homogeneous vector $v\in V$, there is a unique $U_q(\lie g)$-intertwiner $\varphi_{\lambda}^v : M_{\lambda} \to M_{\mu}\otimes V$ such that $<v_{\mu}^*,\varphi_{\lambda}^v v_{\lambda}> = v$ where $\mu = \lambda - \text{weight}(v)$. Given two homogeneous vectors $v, w$ we denote by $\varphi^{v\otimes w}_{\lambda}$ the unique intertwiner given by the composition $(\varphi_{\lambda-\text{wt}(w)}^v\otimes 1)\varphi_{\lambda}^w$. The operator $J(\lambda)(v\otimes w) = <\varphi^{v\otimes w}_{\lambda}>$ is called the fusion matrix. We will denote $\tilde{\Phi}_k^{\varphi_{\lambda}^v}(z)$ by $\tilde{\Phi}_{\lambda,k}^v(z)$ and by $\widetilde{\Phi}_{\lambda,k}^{v\otimes w}(z_1/z_2)$ the  the composition
$$(\tilde{\Phi}_{\lambda-\text{wt}(w),k}^v(z_2)\otimes 1)\circ \tilde{\Phi}_{\lambda,k}^w (z_1): M_{\lambda,k} \to M_{\mu,k}\hat{\otimes} z_2^{-\Delta_2}V[z_2,z_2^{-1}]\hat{\otimes}z_1^{-\Delta_1}V[z_1,z_1^{-1}]$$
The correlation function associated to $\widetilde{\Phi}_{\lambda,k}^{v\otimes w}(z_1/z_2)$ is
$$\mathbf{\Psi}_{\lambda,k}^{v\otimes w}(z_1,z_2) = <v_{\mu}^*,\widetilde{\Phi}_{\lambda,k}^{v\otimes w}(z_1/z_2)v_{\lambda}>$$
It is of the form
$$\mathbf{\Psi}_{\lambda,k}^{v\otimes w}(z_1,z_2) = z_1^{-\Delta_1}z_2^{-\Delta_2}J_{\lambda,k}^{v\otimes w}(z_1/z_2)$$
where $J_{\lambda,k}^{v\otimes w}(z)$ is a meromorphic function on \bb C regular at 0. Then we can consider the fusion matrix with spectral parameter, $J_k(u,\lambda) \in$ End$(V\otimes V)$, given by 
\begin{equation}
J_k(u,\lambda)(v\otimes w) = J_{\lambda,k}^{v\otimes w}(z_1/z_2)
\end{equation}
where $e^{2\pi iu} = z_1/z_2$, and define the exchange matrix
\begin{equation}
R_k(u,\lambda) = J_k(u,\lambda)^{-1}\cal R^{21}|_{V(z_2)\otimes V(z_1)}J_k^{21}(u,\lambda)
\end{equation}

The fact that $R_k$ is really a function of $u$ follows from that it is true for $\cal R^{21}|_{V(z_2)\otimes V(z_1)}$ \cite{efk01}. Then one shows \cite{ev99} that $R_k(u,\lambda)$ is a solution of the quantum dynamical Yang-Baxter equation \eqref{eq:qdybs} (with step 1).


\subsection{Statement of the Main Theorems}

We state now the formula for the exchange matrix $R_k(u,\lambda)$ for the standard representation of $U_q(\lie{sl}_{n+1})$ and give the sequence of gauge transformations that transforms it in Felder's elliptic solution for the proper parameters.

\begin{thm}\label{thm:formula}
$$R_k(u,\lambda) = \,\chi(u,\tau,\gamma)\, \times$$
$$\Big(\sum_{m=0}^{n} E_{m,m}\otimes E_{m,m} + \sum_{l\neq m} \sigma_{m,l}(\lambda,k)\alpha(u,\gamma(\lambda+\rho)_{m,l}) E_{m,m}\otimes E_{l,l} + \beta(u,\gamma(\lambda+\rho)_{m,l}) E_{l,m}\otimes E_{m,l}\Big)$$
where
$$\chi(u,\tau,\gamma) = q^{\frac{n}{n+1}} \frac{\Gamma_e(-u+\tau,(n+1)\gamma,\tau)} {\Gamma_e(-u+\tau+\gamma,(n+1)\gamma,\tau)} \frac{ \Gamma_e(u+\gamma,(n+1)\gamma,\tau)}{\Gamma_e(u,(n+1)\gamma,\tau)}$$
$$\alpha(u,\lambda) = \frac{\vartheta_1(\lambda+\gamma;\tau)}{\vartheta_1(\lambda;\tau)} \frac{\vartheta_1(u;\tau)}{\vartheta_1(u-\gamma;\tau)}\qquad\qquad \beta(u,\lambda) = \frac{\vartheta_1(\gamma;\tau)}{\vartheta_1(\lambda;\tau)} \frac{\vartheta_1(u-\lambda;\tau)}{\vartheta_1(u-\gamma;\tau)}$$
$$\sigma_{l,m}(\lambda,k) = q\frac{\Gamma_p(1+\frac{1}{\kappa}(\lambda+\rho)_{l,m}+\frac{1}{\kappa})} {\Gamma_p(1+\frac{1}{\kappa}(\lambda+\rho)_{l,m})} \frac{\Gamma_p(-\frac{1}{\kappa}(\lambda+\rho)_{l,m})} {\Gamma_p(-\frac{1}{\kappa}(\lambda+\rho)_{l,m}+\frac{1}{\kappa})} \quad\text{and}\quad \sigma_{m,l}(\lambda,k) = \frac{1}{\sigma_{l,m}(\lambda,k)}$$
for $m<l$. Here $p=q^{-2\kappa}$, $|q|<1$, $|p|<1$, $\tau = \frac{\log p}{2\pi i} = -\kappa\frac{\log q}{\pi i}$, $\gamma = -\frac{\tau}{\kappa} = \frac{\log q}{\pi i}$,  $\vartheta_1$ is the standard first theta function \eqref{eq:ftheta}, $\Gamma_q$ is the $q$-gamma function \eqref{eq:thetagamma} and $\Gamma_e$ is the elliptic gamma function \eqref{eq:egamma}.
\end{thm}

\obs We will assume throughout the paper that $|q|<1$ and $|p|<1$. One can manage to write similar formulas for the other three regions $|q|<1<|p|$, $|p|<1<|q|$ and $|q|>1,|p|>1$, writing the infinite products and series involved in the right convergence region, but should be warned that it is not just a matter of changing signs of $\tau$ and $\gamma$ accordingly in theorem \ref{thm:formula}.

\begin{thm}\label{thm:equivalence}
$R_k(u,\lambda)$ is gauge equivalent to Felder's elliptic solution $R_{\gamma,\tau}^{ell}(u,\lambda)$ for $\gamma$ and $\tau$ given in theorem \ref{thm:formula}. The sequence of gauge transformations is
\begin{enumerate}
\item gauge transformation \eqref{eq:g4} with $c(u) = \chi(u,\tau,\gamma)^{-1}$;

\item gauge transformation \eqref{eq:g5} with $a=1$, $b=1/\gamma$ and $\mu = -\rho$;

\item gauge transformation \eqref{eq:g2} with $\varphi_{m,l}(\lambda) = \sigma_{l,m}(\lambda/\gamma-\rho,k)$.
\end{enumerate}
\end{thm}

\obs Etingof and Varchenko \cite{ev97} classified the dynamical R-matrices of $\lie{gl}_{n+1}$ type up to gauge transformations, when $\gamma$ is regarded as a formal parameter (in the case of elliptic quasiclassical limit). The answer is that all solutions are gauge equivalent to Felder's solution. It is believed that in the analytic, rather than formal, situation the classification still holds, which shows why Theorem \ref{thm:equivalence} (with some form of gauge transformations) is to be expected a priori. Nevertheless, these arguments cannot be used for proving Theorem \ref{thm:equivalence}, since the results of \cite{ev97} are proved only in the formal case. Besides, it is useful to have the precise formula for $R_k(u,\lambda)$ and the corresponding sequence of gauge transformations, rather than just an existence result for such gauge transformations. 


\section{The q-\kz Equations}

\subsection{Action of the Universal R-matrix of $U_q(\widetilde{\lie{sl}}_{n+1})$ on Evaluation Representations}\label{sec:er}

For $z \in \bb C^*$ let $D_z$ be the automorphism of $U_q(\hlie g)$ given by $D_z(e_0) = ze_0, D_z(f_0)=z^{-1}f_0$ and the identity on the other generators. Jimbo \cite{jim02} (see also \cite{efk01}) proved the existence of an algebra morphism $p: U_q(\hlie g) \to U_q(\lie g )$ (not a Hopf algebra morphism and only for $\lie g = \lie{sl}_{n+1}$). Then, if $V$ is a finite dimensional representation of $U_q(\lie g)$, the composition $p_z = p\circ D_z$ induces a representation of $U_q(\hlie g)$ in $V$, which is the $q$-analogue of the evaluation representation of $\hlie g$ in $V$ and is denoted by $V(z)$. The element $c$ acts as $0$ on such evaluation representations. Therefore the action of $\cal R$ is well defined on $V(x)\otimes W(y)$ and is given by a power series $\mathbf{R}_{V,W}(x/y)$. It can be shown that this series is convergent in the region $|x| << |y|$, is regular at $0$ and admits a meromorphic continuation to $x/y \in \bb C$. It is convenient to work with the renormalized R-matrix given by the following proposition \cite{efk01} :

\begin{prop}\label{prop:factor}
Let $V,W$ be irreducible evaluation representations of $U_q(\hat{\mathfrak g})$. Then
$$\mathbf{R}_{V,W}(z) = f_{V,W}(z)R_{V,W}(z)$$
where $f_{V,W}$ is a scalar function meromorphic in \bb C, regular at $0$ with $f_{V,W}(0) \neq 0$ and the matrix elements of $R_{V,W}(z)$ are rational functions of $z$ regular at $0$ and such that $R_{V,W}(z)(v_0\otimes w_0) = v_0\otimes w_0$ where $v_0, w_0$ are the highest weight vectors of $V$ and $W$ as $U_q(\lie g)$-modules. Furthermore, if $|q|<1$, $f_{V,W}$ can be represented as
$$f_{V,W}(z) = q^{<\lambda,\mu>}\prod_{j=0}^{\infty} \varrho_{V,W}(q^{2jh\check{}}z)$$
where $\lambda$ and $\mu$ are the highest weights of $V$ and $W$ and $\varrho_{V,W}$ is a rational function such that $\varrho_{V,W}(0)=1$.
\end{prop}

\obs We refer to \cite{efk01} for the proof. We just recall that $\varrho_{V,W}$ is given by the equation 
\begin{equation}\label{eq:ro}
\varrho(z)\,(((R(z)^{-1})^{t_1})^{-1})^{t_1} = (q^{2\rho}\otimes 1)R(q^{2\rho}z)(q^{-2\rho}\otimes 1)
\end{equation}
where $T^{t_1} = \sum a_i^*\otimes b_i$ if $T = \sum a_i\otimes b_i$ is a linear map $V\otimes W \to V\otimes W$.

\begin{prop}\label{prop:unit}\cite{efk01}
Let $U,V,W$ be irreducible evaluation representations of $U_q(\hlie g)$, then $R$ satisfies the quantum Yang-Baxter equation with spectral parameter
\begin{equation}
R_{U,V}(z)R_{U,W}(zw)R_{V,W}(w) = R_{V,W}(w)R_{U,W}(zw)R_{U,V}(z)
\end{equation}
and the unitary property
\begin{equation}\label{eq:unit}
PR_{W,V}(z)PR_{V,W}(z^{-1}) = 1
\end{equation}
where $P$ is the flip map.
\end{prop}

The map $PR_{V,W}(x/y)$ is an intertwiner $V(x)\otimes W(y) \to W(y)\otimes V(x)$. We now calculate the action of this intertwiner when $V=W$ is the standard representation of $U_q(\lie{sl}_{n+1})$. First note that the weight decomposition of $V\otimes V$ is $$V\otimes V = \sum_{m=0}^n V_{m,m} \oplus \sum_{0=m<l}^{n-1} V_{m,l}$$ where $V_{m,m} = \bb C(v_m\otimes v_m)$ and $V_{m,l} = \bb C(v_m\otimes v_l) \oplus \bb C(v_l\otimes v_m)$. $PR(x/y)$ is completely determined by the properties of being an intertwiner and $PR(x/y)(v_0\otimes v_0) = v_0\otimes v_0$.

\begin{prop}\label{prop:urma}
Let $z=x/y$ and $m<l$. The action of $R(z)$ in $V(x)\otimes V(y)$ is given by :
\begin{align*}
R(z)&(v_m\otimes v_m) =  v_m\otimes v_m\\
R(z)&(v_m\otimes v_l) = \xi(z) \,(v_m\otimes v_l) + \eta(z)\, z \,(v_l\otimes v_m)\\
R(z)&(v_l\otimes v_m) = \eta(z) \,(v_m\otimes v_l) + \xi(z)\, (v_l\otimes v_m)
\end{align*}
where
\begin{equation}\label{eq:xieta}
\xi(z) = \frac{1-z}{q-q^{-1}z} \quad\text{and}\quad \eta(z) = \frac{q-q^{-1}}{q-q^{-1}z}
\end{equation}
\end{prop}

The proof is in appendix \ref{ap:urma}. 

It is convenient to write $R(z)$ in the form
\begin{equation*}
R(z) = \sum_{i=0}^n E_{i,i}\otimes E_{i,i} + \xi(z)\Big(\sum_{j>i} E_{i,i}\otimes E_{j,j} + E_{j,j}\otimes E_{i,i}\Big) + \eta(z)\Big( \sum_{j>i} zE_{j,i}\otimes E_{i,j} + E_{i,j}\otimes E_{j,i} \Big)
\end{equation*}

We now calculate $\varrho_{n+1} = \varrho_{V,V}$. We see that $\det R(z) = \xi(z)^2 - z\eta(z)^2$ is the determinant of $R(z)$ in the weight spaces $V_{m,l}$. Then we have
$$(R(z)^{-1})^{t_1} = \sum_{i=0}^n E_{i,i}\otimes E_{i,i} \quad+\quad \frac{\xi(z)}{\det R(z)} \Big(\sum_{j>i} E_{i,i}\otimes E_{j,j} + E_{j,j}\otimes E_{i,i}\Big)\quad$$
$$- \quad \frac{\eta(z)}{\det R(z)}\Big( \sum_{j>i} zE_{i,j}\otimes E_{i,j} + E_{j,i}\otimes E_{j,i} \Big)$$

Now observe that the invariant subspaces of this operator are  $V_= = \oplus \bb C (v_i\otimes v_i)$ and $V_{i\neq j} = \bb C(v_i\otimes v_j)$ and that the action on $V_=$ is given by the matrix
\begin{equation}\label{eq:p}
T_{n+1} =
\begin{pmatrix}
1      & zp(z)  & \cdots & zp(z) \\
p(z)   & 1      & \ddots & \vdots \\
\vdots & \ddots & \ddots & zp(z) \\
p(z)   & \cdots & p(z)   & 1
\end{pmatrix}_{(n+1)\times (n+1)}
\text{where}\qquad
p(z) = - \frac{\eta(z)}{\det R(z)} = \frac{1-q^2}{1-q^2z}
\end{equation}

From here it is easy to see that (apply \eqref{eq:ro} to $v_0\otimes v_0$)
$$\varrho_{n+1} = \frac{\det T_{n+1}}{\det T_n}$$

After calculating this expression (see appendix \ref{ap:fro} for details) we get
\begin{equation}\label{eq:fro}
\varrho_{n+1}(z) = \frac{(1-z)(1-zq^{2(n+1)})}{(1-zq^2)(1-zq^{2n})}
\end{equation} 

Finally we have the expression for $f_{V,V} = f_{n+1}$
\begin{equation}\label{eq:f}
f_{n+1}(z) = q^{\frac{n}{n+1}}\prod_{j=0}^{\infty} \varrho_{n+1}(q^{2j(n+1)}z)
\end{equation}

\obs Recall that $\frac{n}{n+1} = <\omega_1,\omega_1>$ (appendix \ref{ap:useful}). For $\lie{sl}_2$ we have $\varrho_2(z) = \frac{(1-z)(1-q^4z)}{(1-zq^2)^2}$. A more general formula for $\lie{sl}_2$ can be found in \cite{efk01}.

To extend the evaluation representation $V(z)$ to the whole $U_q(\tlie g)$ we consider $z^{-\Delta}V[z,z^{-1}] = V \otimes z^{-\Delta}\bb C[z,z^{-1}]$ where $U_q(\hlie g)$ acts as on $V(z)$ (considering $z$ as a variable and not a number) and $q^d$ will act as $z\frac{d}{dz}$. 


\subsection{Correlation Functions}\label{ssec:cf}

If we begin with $v_i \otimes v_j \in V_{m,l}$ ($i,j \in\{m,l\}$), we have an intertwiner
$$\Phi_{\lambda,k}^{v_i\otimes v_j}(y/x) : M_{\lambda,k} \to M_{\mu,k}\hat{\otimes} V[x,x^{-1}]\hat{\otimes} V[y,y^{-1}]$$
where $\mu = \lambda-(\mu_m+\mu_l)$. The correlation function $\mathbf{\Psi}_{\lambda,k}^{v_i\otimes v_j}(x,y)$ (associated to $\widetilde{\Phi}_{\lambda,k}^{v_i\otimes v_j}(y/x)$) lies in $x^{-\Delta_1}y^{-\Delta_2}V_{m,l}[[y/x]]$ (see \cite{fr01,efk01}), where 
$$\Delta_1= \Delta_k(\lambda-\mu_j) - \Delta_k(\lambda-(\mu_l+\mu_m)) \qquad \Delta_2 = \Delta_k(\lambda) - \Delta_k(\lambda-\mu_j)$$
and the constant term in $V_{m,l}[[y/x]]$ is equal to $J(\lambda)(v_i\otimes v_j)$. 

\obs The fusion matrix $J(\lambda)$ is given in appendix \ref{ap:fus}.
         
The following fundamental theorem was derived in \cite{fr01} (see also \cite{efk01}) :

\begin{thm}\label{thm:qkz}
The correlation functions $\mathbf{\Psi}_{\lambda,k}^{v_i\otimes v_j}$ satisfies the q-\kz equations
$$\mathbf{\Psi}(px,y)= q^{\lambda+\mu+2\rho}_{(1)}  \mathbf{R}^{2,1}(y/x)^{-1}\mathbf{\Psi}(x,y)$$
$$\mathbf{\Psi}(x,py)= \mathbf{R}^{2,1}(py/x) q^{\lambda+\mu+2\rho}_{(2)} \mathbf{\Psi}(x,y)$$
where $p=q^{-2\kappa}$ and $q^{\nu}_{(i)}$ means that the action is on the $i$-th component.
\end{thm}

\obs Recall that we will assume later (to deduce equation \eqref{eq:h} bellow for example) that $|p|<1$, which is really a condition over $k$ in the end. 

Since $\mathbf{\Psi}_{\lambda,k}^{v_i\otimes v_j}$ satisfies the q-\kz equations, we call it a fusion solution.                                                                     

The nontrivial information about the solutions of the q-\kz equations is concentrated in the modified q-\kz equations 
\begin{align}\notag
& \Psi(px,y)= q^{\lambda+\mu+2\rho}_{(1)}  R(x/y)\Psi(x,y)\\\label{eq:mkz}
& \Psi(x,py)= R^{2,1}(py/x) q^{\lambda+\mu+2\rho}_{(2)} \Psi(x,y)
\end{align}
obtained from the original replacing $\mathbf R$ by $R$ (notice we used the unitary property of $R$ to write the modified system). 

Denote by $\mathbf{\Psi}_{\lambda,k}^{m,l}(x,y)$ the solutions of the q-\kz equations in the weight spaces $V_{m,l}$. The solutions of the original equations are given by
\begin{equation}\label{eq:or-m}
\mathbf{\Psi}_{\lambda,k}^{m,l}(x,y) = G(x,y) \Psi_{\lambda,k}^{m,l}(x,y)
\end{equation}
where $\Psi_{\lambda,k}^{m,l}(x,y)$ are the solutions of the modified equations and $G$ is a scalar function satisfying
\begin{equation}\label{eq:eqG}
G(x,y) = f(y/x) G(px,y) \qquad G(x,py) = f(py/x) G(x,y)
\end{equation}
where $f$ is given by \eqref{eq:f}
\begin{equation}\label{eq:g}
f_{n+1}(z) = q^{\frac{n}{n+1}}g(z) \qquad\text{and}\qquad g(z) = \prod_{j=0}^{\infty} \varrho_{n+1}(q^{2j(n+1)}z)
\end{equation}
Then equation \eqref{eq:eqG} can be writen as
$$G(x,y) = p^{-\frac{n}{2\kappa(n+1)}}\,g(y/x) G(px,y) \qquad G(x,py) = p^{-\frac{n}{2\kappa(n+1)}}\,g(py/x) G(x,y)$$
whose solution is
\begin{equation}\label{eq:G}
G(x,y) = (x/y)^{\frac{n}{2\kappa(n+1)}}h(y/x)
\end{equation}
where
\begin{equation}\label{eq:h}
h(z) =\prod_{l=0}^{\infty}g(p^{l+1}\,z)^{-1} = \prod_{j,l=0}^{\infty} \varrho_{n+1}(q^{2j(n+1)}p^{l+1}z)^{-1}
\end{equation}

\obs If we were supposing that $|p|>1$ (but still $|q|<1$), the correct expression for $h(z)$ would be $h(z) =\prod_{l=0}^{\infty}g(p^{-l}\,z)$. In section \ref{sec:eRm} we express $f(z^{-1})\frac{h(z)}{h(z^{-1})}$ in term of elliptic gamma functions.


\subsection{Fusion Solutions of the q-\kz Equations}

We now solve the q-\kz equations for the standard representation of $U_q(\lie{sl}_{n+1})$ taking values in the weight spaces $V_{m,l}, (m\leq l)$. When it does not cause confusion, we will drop the indices $\lambda, k, m, l$. We define
\begin{equation}\label{eq:deltas}
\begin{matrix}
\Delta & = & \Delta_1+\Delta_2 & = & \frac{<2\lambda-(\mu_m+\mu_l)+2\rho,\mu_m+\mu_l>}{2\kappa} \\
\\
\varpi & = & \frac{<\lambda+\rho,\mu_l-\mu_m>}{2\kappa} & = & \frac{1}{2\kappa}(\lambda+\rho)_{l,m}
\end{matrix}
\end{equation}

If $m=l$ one easily sees that the solution of \eqref{eq:mkz} in the space $V_{m,m}$ is
$$\Psi(x,y) = (xy)^{-\Delta/2}v_m\otimes v_m$$
up to multiplication by a pseudo-constant $C(z)$, i.e., $C(pz) = C(z)$. Then $\mathbf{\Psi}(x,y)=G(x,y)\Psi(x,y)$ is the fusion solution in this space since
\begin{equation}\label{eq:1dims}
\mathbf{\Psi}(x,y) = (xy)^{-\Delta/2}(x/y)^{\frac{n}{2\kappa(n+1)}}\, h(y/x)(v_m\otimes v_m)
\end{equation}
and $J(\lambda)(v_m\otimes v_m) = v_m\otimes v_m$.

Now we fix $m<l$. Combining the 2 modified q-\kz equations one easily sees that 
$$\Psi(px,py) = p^{-\Delta}\Psi(x,y)$$
Hence 
$$\Psi(x,y) = (xy)^{-\Delta/2}\psi(x/y)$$
and $\psi(z)$ satisfies
$$\psi(pz)= p^{\Delta/2} p^{-\frac{2\lambda-(\mu_m+\mu_l)+2\rho}{2\kappa}}_{(1)}  R(z)\psi(z)$$
We write 
$$\psi(z) = \psi_1(z)(v_m\otimes v_l) + \psi_2(z)(v_l\otimes v_m)$$
and use proposition \ref{prop:urma} to get the following system
\begin{equation}
\begin{pmatrix}
\psi_1(pz)\\ \psi_2(pz)
\end{pmatrix}
=
\begin{pmatrix}
\xi(z)p^{\varpi}   & \eta(z)p^{\varpi}\\
z\eta(z)p^{-\varpi} & \xi(z)p^{-\varpi}
\end{pmatrix}
\begin{pmatrix}
\psi_1(z)\\ \psi_2(z)
\end{pmatrix}
\end{equation}
where $\xi, \eta$ are given by \eqref{eq:xieta}.

Then one reduces this system to the following equation for $\psi_1$
\begin{equation} \label{eq:2ndord}
(q^{-1}pz-q)\psi_1(p^2z) + \big( p^{\varpi}+p^{-\varpi} - (p^{\varpi+1}+p^{-\varpi})z\big) \psi_1(pz) + (qz-q^{-1}) \psi_1(z) = 0
\end{equation}
which is a second order difference equation as considered in appendix \ref{ap:qhf} (see proposition \ref{prop:qhg}). In our case 
$$A_0 = q^{-1}p, \quad A_1 = -(p^{\varpi+1}+p^{-\varpi}), \quad A_2 = q, \quad B_0 = -q, \quad B_1 = (p^{\varpi}+p^{-\varpi}), \quad B_2 = -q^{-1}$$
and we have two sets of solutions, $(u_i,r_i,s_i,t_i)$, all given up to the sum of an integer multiple of $\frac{2\pi i}{\log p}$ and also up to interchange $r$ and $s$ in each set

\begin{equation}\label{eq:urst}
\begin{matrix}
u_1 = -\varpi+\frac{1}{2\kappa} & r_1 = \frac{1}{\kappa}+1 & s_1 = \frac{1}{\kappa}-2\varpi & t_1 = -2\varpi+1\\
\\
u_2 = \varpi+\frac{1}{2\kappa} & r_2 = \frac{1}{\kappa} & s_2 = \frac{1}{\kappa}+2\varpi+1 & t_2 = 2\varpi+1
\end{matrix}
\end{equation}

Since fusion solutions lie in $x^{-\Delta_1}y^{-\Delta_2}\,V_{m,l}[[y/x]]$, instead of taking the solutions of the corresponding q-hypergeometric equations which are regular at $z=0$, we take their quasimeromorphic solutions given by proposition \ref{prop:infty}. Now we are ready to state
\pagebreak
\begin{prop}\label{prop:sqkz}
The fusion solutions of the q-\kz equations for the standard representation of $U_q(\lie{sl}_{n+1})$ are 
$$\mathbf{\Psi}_{\lambda,k}^{v_m\otimes v_m}(x,y) = (xy)^{-\Delta/2}(x/y)^{\frac{n}{2\kappa(n+1)}}\, h(y/x)(v_m\otimes v_m)$$
if $m=l$ and, for $m<l$
\begin{align*}
& \mathbf{\Psi}_{\lambda,k}^{v_m\otimes v_l}(x,y) = (xy)^{-\Delta/2}(x/y)^{\frac{n}{2\kappa(n+1)}}\, h(y/x) \Big(\psi_1^{(1)}(x/y)(v_m\otimes v_l) + \psi_2^{(1)}(x/y)(v_l\otimes v_m) \Big)\\
& \mathbf{\Psi}_{\lambda,k}^{v_l\otimes v_m}(x,y) = (xy)^{-\Delta/2}(x/y)^{\frac{n}{2\kappa(n+1)}}\, h(y/x) \Big(\psi_1^{(2)}(x/y)(v_m\otimes v_l) + \psi_2^{(2)}(x/y)(v_l\otimes v_m) \Big)
\end{align*}
where
\begin{align*}
& \psi_1^{(1)}(x/y)  =  (x/y)^{\varpi-1/2\kappa}   \,_2\phi_1(p^{1/\kappa}, p^{-2\varpi+1/\kappa}, p^{-2\varpi};p, p^{1-1/\kappa}y/x)\\
& \psi_1^{(2)}(x/y)  =  \epsilon\,(x/y)^{-\varpi-1/2\kappa}  (y/x)\,_2\phi_1(p^{1+1/\kappa},p^{1+2\varpi+1/\kappa},  p^{2(\varpi+1)};p, p^{1-1/\kappa}y/x)\\
& \psi_2^{(i)}(z) =   \frac{p^{-\varpi}(q-q^{-1}z)\psi_1^{(i)}(pz)-(1-z)\psi_1^{(i)}(z)}{q-q^{-1}}
\end{align*}
$$\Delta = \frac{<2\lambda-(\mu_m+\mu_l)+2\rho,\mu_m+\mu_l>}{2\kappa} \qquad
\varpi = \frac{<\lambda+\rho,\mu_l-\mu_m>}{2\kappa} \qquad \epsilon = \frac{q-q^{-1}}{1-p^{-(2\varpi+1)}}$$
$p=q^{-2\kappa}$ and $h(y/x)$ is given by \eqref{eq:h}.
\end{prop}

\begin{proof}
We already know that these are solutions of the q-\kz equations. Hence we are left to check that they are fusion solutions. Begin with the trivial observation that $x^{-\Delta_1}y^{-\Delta_2} = (xy)^{-\frac{\Delta}{2}}(x/y)^{-\frac{\Delta_1-\Delta_2}{2}}$. Now notice that $-\frac{\Delta_1-\Delta_2}{2} = \pm\varpi+\frac{<\mu_m,\mu_l>}{2\kappa}$ respectively for $\mathbf{\Psi}_{\lambda,k}^{v_m\otimes v_l}$ and $\mathbf{\Psi}_{\lambda,k}^{v_l\otimes v_m}$. Recalling the expressions for $<\mu_m,\mu_l>$ in appendix \ref{ap:useful} we are left to check if the degree zero components in $V_{m,l}[[y/x]]$ coincide with the corresponding expectation values. For $m=l$ we have already done it. For the two-dimensional weight spaces, one first notices that, though the expressions for $\psi_2^{(i)}$ seem to contain a term in $x/y$, they actually cancel in both cases. In the case of $\mathbf{\Psi}_{\lambda,k}^{v_l\otimes v_m}$ it is clear that the degree zero coefficient in the direction of $v_m\otimes v_l$ vanishes. Therefore we just have to multiply by a constant ($1$ and $\epsilon$ respectively) to get the predicted expectation values.
\end{proof}


\subsection{Intertwined Fusion Solutions}

Let $\mathbf{\Psi}$ be the correlation function associated to an intertwiner $$(\widetilde{\Phi}_1\otimes1)\widetilde{\Phi}_2:M_{\lambda,k} \to M_{\mu,k}\hat{\otimes}x^{-\Delta_1}V[x,x^{-1}]\hat{\otimes}y^{-\Delta_2}V[y,y^{-1}]$$ and set
$$\phi(y,x) = PR(x/y)\mathbf{\Psi}(x,y)$$
As an immediate consequence of the q-\kz equations for $\mathbf{\Psi}$ we get 
\begin{prop}
$$\phi(py,x) = f(py/x) q^{(2\lambda-(\mu_m+\mu_l)+2\rho)}_{(1)} R(y/x) \phi(y,x)$$
$$\phi(y,px) = \frac{1}{f(y/x)} R^{21}(px/y)q^{(2\lambda-(\mu_m+\mu_l)+2\rho)}_{(2)}  \phi(y,x)$$
where $f$ is given by \eqref{eq:f}.
\end{prop}

One easily recognizes now that $\phi(y,x) = G(x,y) \check{\Psi}(y,x)$, where $G$ is the same given by \eqref{eq:G} and $\check{\Psi}(y,x)$ is the corresponding solution of the modified q-\kz equations intertwined, i.e., with $x,y$ interchanged. Therefore we can calculate the (intertwined) fusion solutions $\check{\mathbf{\Psi}}$ as linear combination (up to multiplication by $\frac{G(y,x)}{G(x,y)}$) of the solutions of the (not intertwined) modified q-\kz equations  using its solutions around $x/y = 0$ after multiplying them by $PR(x/y)$. We recall what these solutions were \eqref{eq:urst}
\begin{align}\notag
&(xy)^{-\Delta/2}\dot{\psi}_1^{(1)}(x/y) =  (xy)^{-\Delta/2}(x/y)^{-\varpi+1/2\kappa}\,_2\phi_1(p^{1+1/\kappa},p^{-2\varpi+1/\kappa}, p^{-2\varpi+1};p, p^{-1/\kappa}x/y)\\\label{eq:as}
&(xy)^{-\Delta/2}\dot{\psi}_1^{(2)}(x/y) =  (xy)^{-\Delta/2}(x/y)^{\varpi+1/2\kappa}\,_2\phi_1(p^{1/\kappa},p^{2\varpi+1+1/\kappa}, p^{2\varpi+1};p, p^{-1/\kappa}x/y)\\
& \dot{\psi}_2^{(i)}(z) =   \frac{p^{-\varpi}(q-q^{-1}z)\dot{\psi}_1^{(i)}(pz)-(1-z)\dot{\psi}_1^{(i)}(z)}{q-q^{-1}}\notag
\end{align}
Now, in the same way we proved proposition \ref{prop:sqkz}, we see that the following are the fusion solutions of the ``intertwined'' q-\kz equations

\begin{equation}\label{eq:i1d}
\check{\mathbf{\Psi}}_{\lambda,k}^{v_m\otimes v_m}(y,x) = (xy)^{-\Delta/2}(y/x)^{\frac{n}{2\kappa(n+1)}}\, h(x/y)(v_m\otimes v_m)
\end{equation}
if $m=l$ and, for $m<l$
\begin{align}\label{eq:i2d1}
& \check{\mathbf{\Psi}}_{\lambda,k}^{v_m\otimes v_l}(y,x) = \check{\epsilon}_1 (xy)^{-\Delta/2}(y/x)^{\frac{n}{2\kappa(n+1)}} h(x/y) \Big(\check{\psi_1}^{(1)}(y/x)(v_m\otimes v_l) + \check{\psi_2}^{(1)}(y/x)(v_l\otimes v_m) \Big)\\\label{eq:i2d2}
& \check{\mathbf{\Psi}}_{\lambda,k}^{v_l\otimes v_m}(y,x) =  \check{\epsilon}_2 (xy)^{-\Delta/2}(y/x)^{\frac{n}{2\kappa(n+1)}}\, h(x/y)   \Big(\check{\psi_1}^{(2)}(y/x)(v_m\otimes v_l) + \check{\psi}_2^{(2)}(y/x)(v_l\otimes v_m) \Big)
\end{align}
where
\begin{align*}
& \check{\psi}_1^{(i)}(y/x)  =  (x/y)\eta(x/y)\dot{\psi}_1^{(i)}(x/y)+\xi(x/y)\dot{\psi}_2^{(i)}(x/y) \\
& \check{\psi}_2^{(i)}(y/x)  =   \xi(x/y)\dot{\psi}_1^{(i)}(x/y)+\eta(x/y)\dot{\psi}_2^{(i)}(x/y)
\end{align*}
$\check{\epsilon}_1 = \frac{q(q-q^{-1})}{p^{-2\varpi}-1}, \check{\epsilon}_2 = q$ and $\xi, \eta$ are given by \eqref{eq:xieta}.


\section{Monodromy and the Elliptic Quantum Dynamical R-Matrix}

\subsection{The Exchange Matrix for Intertwiner Operators}

The $U_q(\hlie g)$-intertwiner $\Phi^{v\otimes w}(x/y):M_{\lambda,k} \to M_{\mu,k}\hat{\otimes}V[y,y^{-1}]\hat{\otimes}V[x,x^{-1}]$, which is analytic in the region $|y|>>|x|$ can be analytically continued to a meromorphic one in the region $|y|<<|x|$. We still denote the analytic continuation by $\Phi^{v\otimes w}(x/y)$. Therefore, the product $PR(y/x)\Phi^{v\otimes w}(x/y)$ is a $U_q(\hlie g)$-intertwiner $M_{\lambda,k} \to M_{\mu,k}\hat{\otimes}V[x,x^{-1}]\hat{\otimes}V[y,y^{-1}]$. Hence it must be of the form $\Phi_{\lambda,k}^{\check{B}(v\otimes w)}(y/x)$ for some operator $\check{B}_k(x/y,\lambda) : V\otimes V \to V\otimes V$. We call the operator
$$B_k(z, \lambda) = \check{B}_k(z,\lambda)P$$
the unitary exchange matrix since it satisfies the unitary condition
\begin{equation}\label{eq:unite}
B_k^{21}(z,\lambda)B_k(z^{-1},\lambda) = 1
\end{equation}                                

The (``non unitary'') exchange matrix constructed in section \ref{sec:fec} is given by
\begin{equation}
R_k(u,\lambda) = f(z^{-1})B_k(z,\lambda)
\end{equation}
where $z =x/y = e^{-2\pi i u}$ and $f$ is given by proposition \ref{prop:factor}. 


\subsection{A Formula for the Exchange Matrix through Monodromy}

We now calculate the exchange matrix $R_k(u,\lambda)$ for the standard representation of $U_q(\lie{sl}_{n+1})$. We begin by the ``twisted'' unitary exchange matrix $\check{B}_k(z,\lambda)$. Since there is no risk of confusion we will drop the sub-indices $\lambda, k$. So we consider the $U_q(\hlie{sl}_{n+1})$-intertwiners $\Phi^{v_m\otimes v_l}(y/x)$ and $\Phi^{\check{B}(v_i\otimes v_j)}(y/x) = PR(y/x)\Phi^{v_i\otimes v_j}(x/y)$ and their expectation values 
$$J^{v_m\otimes v_l}(y/x) = <\Phi^{v_m\otimes v_l}(y/x)> \qquad \check{J}^{v_i\otimes v_j}(y/x) = <\Phi^{\check{B}(v_m\otimes v_l)}(y/x)>$$
From proposition \ref{prop:sqkz} and equations \eqref{eq:as},\eqref{eq:i1d},\eqref{eq:i2d1} and \eqref{eq:i2d2} we find the expectation values in $V_{m,l}$ (we just have to erase the $x^{-\Delta_i} y^{-\Delta_j}$  from the corresponding solutions of the q-\kz equations). For the 1-dimensional weight spaces $V_{m,m}$ we have 
$$J(y/x) = h(y/x)v_m\otimes v_m \qquad\text{while}\qquad \check{J}(y/x) = h(x/y)v_m\otimes v_m$$
Hence
$$\check{B}_k(x/y,\lambda)(v_m\otimes v_m) = \frac{h(x/y)}{h(y/x)}\,v_m\otimes v_m$$
Now, for $m<l$
$$J^{v_m\otimes v_l}(y/x) = h(y/x)\big(J_1^{(1)}(x/y)(v_m\otimes v_l) + J_2^{(1)}(x/y)(v_l\otimes v_m)\big)$$
$$J^{v_l\otimes v_m}(y/x) = h(y/x)\big(J_1^{(2)}(x/y)(v_m\otimes v_l) + J_2^{(2)}(x/y)(v_l\otimes v_m)\big)$$
where
\begin{align*}
& J_1^{(1)}(z) = \, _2\phi_1(p^{1/\kappa}, p^{-2\varpi+1/\kappa}, p^{-2\varpi};p, p^{1-1/\kappa}z^{-1}) \\
& J_1^{(2)} (z)= \epsilon\,z^{-1}\, _2\phi_1(p^{1+1/\kappa}, p^{1+2\varpi+1/\kappa}, p^{2(\varpi+1)};p, p^{1-1/\kappa}z^{-1})\\
& J_2^{(1)}(z) = \frac{p^{-1/2\kappa}(q-q^{-1}z)J_1^{(1)}(pz)-(1-z)J_1^{(1)}(z)}{q-q^{-1}}
\\
& J_2^{(1)}(z) = \frac{p^{-2\varpi-1/2\kappa}(q-q^{-1}z)J_1^{(2)}(pz)-(1-z)J_1^{(2)}(z)}{q-q^{-1}}
\end{align*}
and
$$\check{J}^{v_m\otimes v_l}(y/x) = h(x/y)\big(\check{J}_1^{(1)}(x/y)(v_m\otimes v_l) + \check{J}_2^{(1)}(x/y)(v_l\otimes v_m)\big)$$
$$\check{J}^{v_l\otimes v_m}(y/x) = h(x/y)\big(\check{J}_1^{(2)}(x/y)(v_m\otimes v_l) + \check{J}_2^{(2)}(x/y)(v_l\otimes v_m)\big)$$
where
\begin{align*}
& \check{J}_1^{(1)}(z) = \check{\epsilon}_1\, _2\phi_1(p^{1+1/\kappa}, p^{-2\varpi+1/\kappa}, p^{-2\varpi+1};p, p^{-1/\kappa}z) \\
& \check{J}_1^{(2)} (z)= \check{\epsilon}_2\, _2\phi_1(p^{1/\kappa}, p^{1+2\varpi+1/\kappa}, p^{2\varpi+1};p, p^{-1/\kappa}z)\\
& \check{J}_2^{(1)}(z) = \frac{p^{-2\varpi+1/2\kappa}(q-q^{-1}z)\check{J}_1^{(1)}(pz)-(1-z)\check{J}_1^{(1)}(z)}{q-q^{-1}}
\\
& \check{J}_2^{(2)}(z) = \frac{p^{1/2\kappa}(q-q^{-1}z)\check{J}_1^{(2)}(pz)-(1-z)\check{J}_1^{(2)}(z)}{q-q^{-1}}
\end{align*}                        

Now we use theorem \ref{thm:mon} to write $\check{J}_1^{(i)}$ as a combination of $J_1^{(i)}$. Of course that the same will automatically be done for $\check{J}_2^{(i)}$.
\begin{align*}
\check{J}_1^{(1)}(z) = &\, \frac{\Gamma_p(1-2\varpi)\Gamma_p(1+2\varpi)}{\Gamma_p(1+1/\kappa)\Gamma_p(1-1/\kappa)} \, \frac{\Theta(zp^{-2\varpi};p)}{\Theta(zp^{-1/\kappa};p)}\,\check{\epsilon}_1 J_1^{(1)}(z)\\
+ &\, \frac{\Gamma_p(1-2\varpi)\Gamma_p(-1-2\varpi)} {\Gamma_p(-2\varpi+1/\kappa)\Gamma_p(-2\varpi-1/\kappa)} \, \frac{\Theta(zp;p)}{\Theta(zp^{-1/\kappa};p)}\,\frac{\check{\epsilon}_1 z}{\epsilon} J_1^{(2)}(z)\\
\\
\check{J}_1^{(2)}(z) = &\, \frac{\Gamma_p(1+2\varpi)\Gamma_p(1+2\varpi)} {\Gamma_p(2\varpi+1+1/\kappa)\Gamma_p(2\varpi+1-1/\kappa)} \, \frac{\Theta(z;p)}{\Theta(zp^{-1/\kappa};p)}\, \check{\epsilon}_2 J_1^{(1)}(z)\\
+ &\, \frac{\Gamma_p(1+2\varpi)\Gamma_p(-1-2\varpi)}{\Gamma_p(1/\kappa)\Gamma_p(-1/\kappa)} \, \frac{\Theta(zp^{1+2\varpi};p)}{\Theta(zp^{-1/\kappa};p)}\, \frac{\check{\epsilon}_2 z}{\epsilon} J_1^{(2)}(z)
\end{align*}

Finally, in the basis $\{v_m\otimes v_l,\, v_l\otimes v_m\}$ for $V_{m,l}$, we have
\begin{gather}\notag
\check{B}_k(z,\lambda) = \frac{h(z)}{h(z^{-1})} \times\\\label{eq:tem}
\begin{pmatrix}
\frac{\Gamma_p(1-2\varpi)\Gamma_p(1+2\varpi)}{\Gamma_p(1+1/\kappa)\Gamma_p(1-1/\kappa)} \, \frac{\Theta(zp^{-2\varpi};p)}{\Theta(zq^2;p)} \,\check{\epsilon}_1 & \frac{\Gamma_p(1+2\varpi)\Gamma_p(1+2\varpi)} {\Gamma_p(2\varpi+1+1/\kappa)\Gamma_p(2\varpi+1-1/\kappa)} \, \frac{\Theta(z;p)}{\Theta(zq^2;p)} \,\check{\epsilon}_2
\\
\\
\frac{\Gamma_p(1-2\varpi)\Gamma_p(-1-2\varpi)} {\Gamma_p(-2\varpi+1/\kappa)\Gamma_p(-2\varpi-1/\kappa)} \, \frac{\Theta(zp;p)}{\Theta(zq^2;p)} \frac{\check{\epsilon}_1 z}{\epsilon} & \frac{\Gamma_p(1+2\varpi)\Gamma_p(-1-2\varpi)}{\Gamma_p(1/\kappa)\Gamma_p(-1/\kappa)} \, \frac{\Theta(zp^{1+2\varpi};p)}{\Theta(zq^2;p)}\frac{\check{\epsilon}_2 z}{\epsilon}
\end{pmatrix}
\end{gather}
with $z=x/y = e^{-2\pi iu}$. Then $R_k(u,\lambda)$ has the form
$$R_k(u,\lambda) = f(z^{-1})\frac{h(z)}{h(z^{-1})} \Big( \sum_{m=0}^{n} E_{m,m}\otimes E_{m,m} + \sum_{m\neq l} \alpha_k^{m,l}(z,\lambda) E_{m,m}\otimes E_{l,l} + \beta_k^{m,l}(z,\lambda) E_{l,m}\otimes E_{m,l} \Big)$$
where, for $m<l$
\begin{align*}
 & \alpha_k^{m,l}(z,\lambda) = \frac{\Gamma_p(1+2\varpi)\Gamma_p(1+2\varpi)} {\Gamma_p(2\varpi+1+1/\kappa)\Gamma_p(2\varpi+1-1/\kappa)} \, \frac{\Theta(z;p)}{\Theta(zq^2;p)} \,\check{\epsilon}_2\\
 & \alpha_k^{l,m}(z,\lambda) = \frac{\Gamma_p(1-2\varpi)\Gamma_p(-1-2\varpi)} {\Gamma_p(-2\varpi+1/\kappa)\Gamma_p(-2\varpi-1/\kappa)} \, \frac{\Theta(zp;p)}{\Theta(zq^2;p)} \frac{\check{\epsilon}_1 z}{\epsilon}\\
 & \beta_k^{m,l}(z,\lambda) = \frac{\Gamma_p(1+2\varpi)\Gamma_p(-1-2\varpi)}{\Gamma_p(1/\kappa)\Gamma_p(-1/\kappa)} \, \frac{\Theta(zp^{1+2\varpi};p)}{\Theta(zq^2;p)}\frac{\check{\epsilon}_2 z}{\epsilon}
\\
 & \beta_k^{l,m}(z,\lambda) = \frac{\Gamma_p(1-2\varpi)\Gamma_p(1+2\varpi)}{\Gamma_p(1+1/\kappa)\Gamma_p(1-1/\kappa)} \, \frac{\Theta(zp^{-2\varpi};p)}{\Theta(zq^2;p)} \,\check{\epsilon}_1\\
\end{align*}

Using identities \eqref{eq:ids}, \eqref{eq:zids}  and recalling the expressions for $\epsilon, \check{\epsilon}_1, \check{\epsilon}_2$ we get expressions for $\beta_k^{i,j}$ written exclusively in terms of $\Theta$ functions
\begin{equation}\label{eq:beta}
\beta_k^{m,l}(z,\lambda) = -p^{-2\varpi} \frac{\Theta(q^2;p)}{\Theta(p^{-2\varpi};p)} \frac{\Theta(zp^{2\varpi};p)}{\Theta(zq^2;p)}
 \qquad\qquad \beta_k^{l,m}(z,\lambda) = \frac{\Theta(q^2;p)}{\Theta(p^{-2\varpi};p)} \frac{\Theta(zp^{-2\varpi};p)}{\Theta(zq^2;p)}
\end{equation}
For $\alpha_k^{i,j}$ we only get
\begin{equation}\label{eq:alpha}
\alpha_k^{m,l}(z,\lambda) = \sigma_{m,l} \, q^2\frac{\Theta(p^{-2\varpi}q^{-2};p)}{\Theta(p^{-2\varpi};p)} \frac{\Theta(z;p)}{\Theta(zq^2;p)}
\quad\quad \alpha_k^{l,m}(z,\lambda) = \sigma_{l,m} \frac{\Theta(p^{-2\varpi}q^{2};p)}{\Theta(p^{-2\varpi};p)} \frac{\Theta(z;p)}{\Theta(zq^2;p)}
\end{equation}
where
\begin{equation}\label{eq:sigma}
\sigma_{l,m}(\lambda,k) = q\frac{\Gamma_p(1+2\varpi+1/\kappa)}{\Gamma_p(1+2\varpi)} \frac{\Gamma_p(-2\varpi)}{\Gamma_p(-2\varpi+1/\kappa)} \quad\text{and}\quad \sigma_{m,l}(\lambda,k) = \frac{1}{\sigma_{l,m}(\lambda,k)}
\end{equation}


\subsection{The Elliptic R-Matrix}\label{sec:eRm}

We now calculate the gauge transformations that show that the exchange matrix for the standard representation of $U_q(\lie{sl}_{n+1})$ is gauge equivalent to Felder's elliptic solution for the the proper parameters.

The exchange matrix is a solution of \eqref{eq:qdybs} with step 1. The first gauge transformation we apply is \eqref{eq:g4}. It just removes the function $f(z^{-1})\frac{h(z)}{h(z^{-1})}$ and doesn't change the step. For the next steps we need to write the exchange matrix in terms of $\vartheta_1$. As usual we set $z=e^{-2\pi i u}$ and also $p=e^{2\pi i\tau}$,i.e., 
$$\tau = \frac{\log p}{2\pi i} = -\frac{\log q}{\pi i}\, \kappa$$
and use \eqref{eq:thetas} and the fact that $\vartheta_1$ is an odd function to obtain, for $m<l$
\begin{align*}
& \alpha_k^{m,l}(z,\lambda) = \sigma_{m,l} \frac{\vartheta_1(2\varpi\tau-\tau/\kappa;\tau)}{\vartheta_1(2\varpi\tau;\tau)} \frac{\vartheta_1(u;\tau)}{\vartheta_1(u+\tau/\kappa;\tau)}\\
& \alpha_k^{l,m}(z,\lambda) = \sigma_{l,m} \frac{\vartheta_1(-2\varpi\tau-\tau/\kappa;\tau)}{\vartheta_1(-2\varpi\tau;\tau)} \frac{\vartheta_1(u;\tau)}{\vartheta_1(u+\tau/\kappa;\tau)} \\
& \beta_k^{m,l}(z,\lambda) =  \frac{\vartheta_1(-\tau/\kappa;\tau)}{\vartheta_1(2\varpi\tau;\tau)} \frac{\vartheta_1(u-2\varpi\tau;\tau)}{\vartheta_1(u+\tau/\kappa;\tau)}
\\
& \beta_k^{l,m}(z,\lambda) =  \frac{\vartheta_1(-\tau/\kappa;\tau)}{\vartheta_1(-2\varpi\tau;\tau)} \frac{\vartheta_1(u+2\varpi\tau;\tau)}{\vartheta_1(u+\tau/\kappa;\tau)} \\
\end{align*}

Then set 
$$\gamma = -\tau/\kappa = \frac{\log q}{\pi i}$$
and apply gauge transformation \eqref{eq:g5} with $a=1$, $b=1/\gamma$ and $\mu = -\rho$. This changes the step to $\gamma$.
One easily checks now that
$$\varphi_{m,l}(\lambda) = \sigma_{l,m}(\lambda/\gamma-\rho,k)$$
is a $\gamma$-closed 2-form and apply gauge transformation \eqref{eq:g2} to finally obtain Felder's solution for these $\gamma$ and $\tau$.

Now that we have identified $\gamma$ and $\tau$, it is interesting to observe how they appear also in the function 
$$\chi(u,\tau,\gamma) = f(z^{-1})\frac{h(z)}{h(z^{-1})}$$
with the usual identifications $z = e^{-2\pi iu}, p = e^{2\pi i\tau} \so q=e^{2\pi i\gamma/2}$. Recalling the expression for $f(z)$ \eqref{eq:f} and $h(z)$ \eqref{eq:h} we see that $\chi$ can be written as
$$\chi(u,\tau,\gamma) = q^{\frac{n}{n+1}}P_1P_2P_3P_4$$
where
\begin{align*}
& P_1 = \prod_{j,l=0}^{\infty} \frac{1-q^{2j(n+1)}p^{l+1}q^{2(n+1)}z^{-1}}{1-q^{2j(n+1)}p^{l+1}z}\, (1-q^{2j(n+1)}q^{2(n+1)}z^{-1})\\
& P_2 = \prod_{j,l=0}^{\infty} \frac{1-q^{2j(n+1)}p^{l+1}z^{-1}}{1-q^{2j(n+1)}p^{l+1}q^{2(n+1)}z}\,(1-q^{2j(n+1)}z^{-1})\\
\end{align*}
\begin{align*}
& P_3 = \prod_{j,l=0}^{\infty} \frac{1-q^{2j(n+1)}p^{l+1}q^2z}{1-q^{2j(n+1)}p^{l+1}q^{2n}z^{-1}} \frac{1}{(1-q^{2j(n+1)}q^{2n}z^{-1})}\\
& P_4 = \prod_{j,l=0}^{\infty} \frac{1-q^{2j(n+1)}p^{l+1}q^{2n}z}{1-q^{2j(n+1)}p^{l+1}q^2z^{-1}} \frac{1}{(1-q^{2j(n+1)}q^{2}z^{-1})}
\end{align*}
Each of these products can be written in terms of a single elliptic gamma function \eqref{eq:egamma} as follows
\begin{align*}
P_1 = &\, \prod_{j,l=0}^{\infty} \frac{1-q^{(j+1)2(n+1)}p^{l}z^{-1}}{1-q^{j2(n+1)}p^{l+1}z} =  \prod_{j,l=0}^{\infty} \frac{1-q^{(j+1)2(n+1)}p^{l+1}(p^{-1}z^{-1})} {1-q^{j2(n+1)}p^l(p^{-1}z^{-1})^{-1}} = \Gamma_e(-u+\tau,(n+1)\gamma, \tau)\\
P_2 = &\, \prod_{j,l=0}^{\infty} \frac{1-q^{j2(n+1)}p^{l}z^{-1}}{1-q^{(j+1)2(n+1)}p^{l+1}z} = \frac{1}{\Gamma_e(u,(n+1)\gamma, \tau)}\\
P_3 = &\, \prod_{j,l=0}^{\infty} \frac{1-q^{j2(n+1)}p^{l+1}q^2z} {1-q^{j2(n+1)}p^lz^{-1}} =  \prod_{j,l=0}^{\infty} \frac{1-q^{j2(n+1)}p^l(pq^2z)} {1-q^{(j+1)2(n+1)}p^{l+1}(pq^2z)^{-1}} \\
= &\,\frac{1}{\Gamma_e(-u+\tau+\gamma,(n+1)\gamma, \tau)}\\
P_4 = &\, \prod_{j,l=0}^{\infty} \frac{1-q^{j2(n+1)}p^{l+1}q^{2n}z} {1-q^{j2(n+1)}p^lq^2z^{-1}} =  \prod_{j,l=0}^{\infty} \frac{1-q^{(j+1)2(n+1)}p^{l+1}(q^2z^{-1})^{-1}} {1-q^{j2(n+1)}p^l(q^2z^{-1})} \\ =&\, \Gamma_e(u+\gamma,(n+1)\gamma, \tau)
\end{align*}

The elliptic gamma function was studied in \cite{felv99}.


\appendix
\section*{Appendix}

\section{{\large Second Order Difference Equations and the q-Hypergeometric Function}}\label{ap:qhf}

Consider a second order difference equation of the form
\begin{equation}
(A_0z+B_0)f(p^2z) + (A_1z+B_1)f(pz) + (A_2z+B_2)f(z) = 0
\end{equation}
For generic values of the coefficients this equation reduces to the Heine's q-hypergeometric equation (see \cite{efk01})
\begin{equation}\label{eq:qhg}
(p^{r+s}z - p^{t-1})f(p^2z) + (-(p^r+p^s)z + p^{t-1} +1)f(pz) + (z-1)f(z) = 0
\end{equation}

\begin{prop}\label{prop:qhg}
Let $r_j, s_j, t_j, u_j, j=1,2$, be the solutions of
$$p^{2u}B_0 + p^uB_1 + B_2 = 0,\quad p^{2u}\frac{A_0}{A_2} = p^{r+s}, \quad p^u\frac{A_1}{A_2} = -(p^r+p^s), \quad p^{2u}\frac{B_0}{B_2} = p^{t-1}$$
Then the functions
$$f^{(j)}(z) = z^{u_j}\, _2\phi_1(p^{r_j}, p^{s_j}, p^{t_j};p, -zA_2/B_2)$$
are the only solutions of \eqref{eq:2ndord} of the form $f(z) = z^ug(z)$ for some $g(z)$ regular in a neighborhood of 0 such that $g(0)=1$, where $_2\phi_1$ is Heine's q-hypergeometric function.
\end{prop}

The function $f(z) = \,_2\phi_1(p^r,p^s,p^t;p,z^{-1})$ satisfies
\begin{equation}
(1-p^2z)f(p^2z) + \big( (p^{t+1}+p^2)z-(p^r+p^s)\big)f(pz) + (p^{r+s}-p^{t+1}z)f(z) = 0
\end{equation}
which is obtained from Heine's q-hypergeometric equation \eqref{eq:qhg} after replacing $z$ by $z^{-1}$. This equation can also in general be reduced to \eqref{eq:qhg} with new parameters $(r,s,t)$ given by
\begin{equation}
(r',s',t') = (r, r-t+1, r-s+1) \qquad\text{or}\qquad (r'',s'',t'') = (s-t+1, s, s-r+1) 
\end{equation}
and one concludes

\begin{prop}\label{prop:infty}
The functions
$$g_1(z) = z^{-r}\, _2\phi_1(p^r,p^{r-t+1},p^{r-s+1};p, p^{t+1-r-s}z^{-1})$$
$$g_2(z) =  z^{-s}\, _2\phi_1(p^{s-t+1},p^s,p^{s-r+1};p, p^{t+1-r-s}z^{-1})$$
are solutions of \eqref{eq:qhg} quasimeromorphic on $\bb C^*\cup \infty$
\end{prop}

Let $\Gamma_q$ be the q-gamma function and $\Theta(z;q)$ be given by \eqref{eq:thetagamma}. Then we have the following theorem for the q-hypergeometric functions \cite{gara01,efk01}

\begin{thm}\label{thm:mon}
$_2\phi_1(p^r,p^s,p^t;p,z) = \Lambda(z)g_1(z) + \Omega(z)g_2(z)$ where
$$\Lambda(z) = \frac{\Gamma_p(t)\Gamma_p(s-r)}{\Gamma_p(s)\Gamma_p(t-r)}\, \frac{\Theta(zp^r;p)}{\Theta(z;p)}z^r$$
$$\Omega(z) = \frac{\Gamma_p(t)\Gamma_p(r-s)}{\Gamma_p(r)\Gamma_p(t-s)}\, \frac{\Theta(zp^s;p)}{\Theta(z;p)}z^s$$
\end{thm}


\section{{\large Theta, Elliptic and q-Gamma Functions}}\label{ap:tgf}

We recall the expressions and identities used
\begin{equation}\label{eq:thetagamma}
\Theta(z;q) = \prod_{j\geq 0} (1-zq^j)(1-z^{-1}q^{j+1})(1-q^{j+1}) \quad\quad\quad \Gamma_q(a) = (1-q)^{1-a}\prod_{j\geq 0} \frac{1-q^{j+1}}{1-q^{j+a}}
\end{equation}
Then
\begin{equation}\label{eq:zids}
z\Theta(qz;q) = -\Theta(z;q) \quad\quad z\Theta(z;q) = -q\Theta(q^{-1}z;q) \quad\quad \Theta(qz;q) = \Theta(z^{-1};q)
\end{equation}

Let $C(q) = (1-q)(q;q)_{\infty}^3 = (1-q)\big(\prod_{j\geq 0} (1-q^{j+1})\big)^3$
\begin{equation}\label{eq:bids}
\Gamma_q(a)\Gamma_q(1-a) = \frac{C(q)}{\Theta(q^a;q)} \quad\quad\quad \Gamma_q(1+a) = \{a\}_q\Gamma_q(a) = \frac{1-q^a}{1-q}\,\Gamma_q(a)
\end{equation}
From this we get
\begin{align}\notag
& \Gamma_q(a)\Gamma_q(2-a) = \frac{\{1-a\}_qC(q)}{\Theta(q^a;q)} = \frac{\{a-1\}_qC(q)}{\Theta(q^{2-a};q)}\\\label{eq:ids}
& \Gamma_q(a)\Gamma_q(-a) = \frac{C(q)}{\{-a\}_q\Theta(q^a;q)} = \frac{C(q)}{\{a\}_q\Theta(q^{-a};q)}
\end{align}

Now the relation between the theta function $\Theta$ and the standard first theta function
\begin{align}\notag
\vartheta_1(u;\tau) & = - \sum_{j=-\infty}^{\infty} e^{\pi i(j+1/2)^2\tau+2\pi i(j+1/2)(u+1/2)}\\\label{eq:ftheta}
 & = -i\,e^{\frac{\pi i\tau}{4}} e^{\pi i u} \prod_{j=0}^{\infty} \big( 1-e^{-2\pi iu}(e^{2\pi i\tau})^{j}\big) \big(1-e^{2\pi i u} (e^{2\pi i\tau})^{j+1} \big) \big( 1 - (e^{2\pi i\tau})^{j+1}\big)
\end{align}
As an immediate consequence of the first expression we see that $\vartheta_1$ is an odd function on $u$. Letting $z=e^{-2\pi iu}$ and $q = e^{2\pi i\tau}$ we get
\begin{equation}\label{eq:thetas}
\Theta(z;q) = iq^{-1/8}z^{1/2} \vartheta_1(u;\tau)
\end{equation}

Finally the elliptic Gamma function
 is given by
\begin{equation}\label{eq:egamma}
\Gamma_e(u,\zeta,\sigma) = \prod_{j,l=0}^{\infty} \frac{1-e^{2\pi i((j+1)\zeta+ (l+1)\sigma -u)}}{1- e^{ 2\pi i(j\zeta+l\sigma+u)}}
\end{equation} 

\section{{\large Proof of Proposition \ref{prop:urma}}}\label{ap:urma}

Using Jimbo's morphism, one sees that the only factor of $p(f_0)$ that acts non trivially on the standard representation is $e_1e_2\dots e_n = E_{0,n}$. For brevity of notation, we will write $f_0$ instead of $p_z(f_0)$. Let $m<l$ and define
$$\begin{matrix}
PR(z)(v_m\otimes v_l) & = & \alpha_{m,l}(z)(v_m\otimes v_l) & + & \beta_{m,l}(z)(v_l\otimes v_m)\\                      
PR(z)(v_l\otimes v_m) & = &\alpha_{l,m}(z)(v_m\otimes v_l)  & + & \beta_{l,m}(z)(v_l\otimes v_m)\\
PR(z)(v_l\otimes v_l) & = & \gamma_l(z)(v_l\otimes v_l)& &
\end{matrix}$$
Fix $l \neq 0,n$. We have
$$PR(x/y)\Delta(f_0)(v_n\otimes v_l) = PR(x/y)(x^{-1}v_0\otimes v_l) = x^{-1}\alpha_{0,l}(x/y)(v_0\otimes v_l) + x^{-1}\beta_{0,l}(x/y)(v_l\otimes v_0)$$
By the intertwining property this is also equal to
$$\begin{matrix}
\Delta(f_0)PR(x/y)(v_n\otimes v_l) & = & \Delta(f_0)\Big(\alpha_{n,l}(x/y)(v_l\otimes v_n) + \beta_{n,l}(x/y)(v_n\otimes v_l)\Big)\\
& &\\
& = & x^{-1}\alpha_{n,l}(x/y)(v_l\otimes v_0) + y^{-1}\beta_{n,l}(x/y)(v_0\otimes v_l)
\end{matrix}$$
and we conclude
\begin{equation}\label{eq:12}
\alpha_{0,l} = z\beta_{n,l} \qquad \beta_{0,l} = \alpha_{n,l}
\end{equation}

The same calculation with $l=0$ will give us
\begin{equation}\label{eq:3}
q\alpha_{n,0} + z\beta_{n,0} = 1
\end{equation}

Now we fix $m<l$ and perform the same calculation with
$$PR(x/y)\Delta(e_m)(v_m\otimes v_l) = \Delta(e_m)PR(x/y)(v_m\otimes v_l)$$
to get
$$\alpha_{m-1,l} = \alpha_{m,l} \qquad \beta_{m-1,l} = \beta_{m,l}$$
Therefore
\begin{equation}\label{eq:56}
\alpha_{m,l} = \alpha_{0,l} \qquad \beta_{m,l} = \beta_{0,l} \qquad m<l
\end{equation}

Also, for $PR(x/y)\Delta(e_l)(v_l\otimes v_l) = \Delta(e_l)PR(x/y)(v_l\otimes v_l)$, we get
\begin{equation}\label{eq:4}
\gamma_l = \alpha_{l-1,l} + q\alpha_{l,l-1} = q^{-1}\beta_{l-1,l} + \beta_{l,l-1}
\end{equation}

For $m>l$ we do the same thing replacing $e_m$ by $f_{m+1}$ and find
\begin{equation}\label{eq:78}
\alpha_{m,l} = \alpha_{n,l} \qquad \beta_{m,l} = \beta_{n,l} \qquad m>l
\end{equation}

Now $PR(x/y)\Delta(e_lf_0)(v_n \otimes v_l) = \Delta(e_lf_0)PR(x/y)$ with $l\neq 0,n$ will give
\begin{equation}\label{eq:9}
q\alpha_{n,1} + z\beta_{n,1} = 1
\end{equation}
\begin{equation}\label{eq:9'}
\alpha_{0,l-1} = z\beta_{n,l} \qquad \alpha_{n,l} = \beta_{0,l-1}\quad l\neq 0,1,n
\end{equation}
Then we combine \eqref{eq:9'}, \eqref{eq:12} and \eqref{eq:56} to get
\begin{equation}\label{eq:1011}
\alpha_{m,l} = \alpha_{0,1} \qquad \beta_{m,l} = \beta_{0,1} \qquad m<l
\end{equation}

Once more the same argument with $PR(x/y)\Delta(f_1)(v_0 \otimes v_0)$ provides
\begin{equation}\label{eq:1213}
q\alpha_{1,0} + \alpha_{0,1}= 1 \qquad q\beta_{1,0} + \beta_{0,1} = q
\end{equation}

Now apply \eqref{eq:78} with $m=1$ in \eqref{eq:3} and use \eqref{eq:1213} to get
$$\alpha_{0,1} = (1-q^{-1}\beta_{0,1})z$$
Then we use in sequence on the right hand side \eqref{eq:12}, \eqref{eq:9} and \eqref{eq:12} to finally obtain
\begin{equation}\label{eq:a01}
\alpha_{0,1}(z) = \frac{q-q^{-1}}{q-q^{-1}z}\, z
\end{equation}
and \eqref{eq:1213} implies
\begin{equation}\label{eq:a10}
\alpha_{1,0}(z) = \frac{1-z}{q-q^{-1}z}
\end{equation}

Going back to the intertwining program we also get
\begin{equation}\label{eq:1617}
\alpha_{m,l} = \alpha_{1,0} \qquad \beta_{m,l} = \beta_{1,0} \qquad m>l
\end{equation}
and we are done with all $\alpha_{i,j}$ and $\beta_{i,j}$ after using \eqref{eq:12}, \eqref{eq:56}, \eqref{eq:78}. Equation \eqref{eq:4} will then tell $\gamma_l = 1$\\
\cqd


\section{{\large The Expression for $\varrho_{n+1}$}}\label{ap:fro}

Recall the definition of the matrix
$$T_{n+1} =
\begin{pmatrix}
1      & zp(z)  & \cdots & zp(z) \\
p(z)   & 1      & \ddots & \vdots \\
\vdots & \ddots & \ddots & zp(z) \\
p(z)   & \cdots & p(z)   & 1
\end{pmatrix}_{(n+1)\times (n+1)}$$

We will calculate $\det T_{n+1}$. First observe that $\det T_{n+1}$ is equal to the determinant of the following matrix
$$\begin{pmatrix}
1      & zp(z)   & \cdots  & \cdots & zp(z)\\
p(z)-1 & 1-zp(z) & 0       & \cdots & 0\\
0      & p(z)-1  & 1-zp(z) & \ddots & \vdots\\
\vdots & \ddots  & \ddots  & \ddots & 0 \\
0      & \cdots  & 0       & p(z)-1 & 1-zp(z)
\end{pmatrix}$$
obtained from $T_{n+1}$ by substituting the i-th row by itself minus the (i-1)-th. We calculate the determinant using the last collum to obtain
$$\det T_{n+1} = (1-zp(z))\det T_n + (-1)^n zp(z)(p(z)-1)^n$$

Similarly (just transpose $T_{n+1}$ and do the same thing) we see that $\det T_{n+1}$ is also equal to the determinant of
$$\begin{pmatrix}
1       & p(z)    & \cdots  & \cdots & p(z)\\
zp(z)-1 & 1-p(z)  & 0       & \cdots & 0\\
0       & zp(z)-1 & 1-p(z)  & \ddots & \vdots\\
\vdots & \ddots   & \ddots  & \ddots & 0 \\
0      & \cdots   & 0       & zp(z)-1 & 1-p(z)
\end{pmatrix}$$
and
$$\det T_{n+1} = (1-p(z))\det T_n + (-1)^n p(z)(zp(z)-1)^n$$
Therefore
$$\det T_n = (-1)^n\, \frac{(zp(z)-1)^n - z\,(p(z)-1)^n}{1-z}$$

For our $p(z) = \frac{1-q^2}{1-q^2z}$ we finally get
$$\det T_n = \frac{(1-q^{2n}z)(1-z)^{n-1}}{(1-q^2z)^n}$$


\section{{\large Useful Formulas}}\label{ap:useful}

Recall that $\lie h$ can be identified with the subspace of $\bb C^{n+1}$ where the sum of the entries of its vectors vanishes. Then we have
\begin{equation}\label{eq:sroots}
h_i \quad=\quad \frac{1}{n+1}(0,0,\dots,1_{i-1},-1_i, 0,\dots,0)
\end{equation}
The invariant symmetric bilinear form is just $< (\lambda_0, \dots, \lambda_n),(\mu_0, \dots, \mu_n)> = \sum \lambda_i\mu_i$. Using $<,>$ to identify $\lie h^*$ with $\lie h$ we find that  the simple roots coincide with the $h_i$ and that the fundamental weights are
\begin{equation}\label{eq:fweights}
\omega_i = \frac{1}{n+1}(n+1-i,\dots, n+1-i, -i,\dots, -i)
\end{equation}
where the change occurs from the $i$-th  to the $i$+1-th coordinate.

Recall that $\{ v_m\}\, (m=0, \dots, n)$ is the canonical basis of $\bb C^{n+1}$ and that $\mu_m$ is the weight of $v_m$. We know that $\mu_0 = \omega_1$ and $\mu_n = -\omega_n$. The expressions for the other $\mu_m$ are
\begin{equation}\label{eq:weights}
\mu_m = \frac{1}{n+1}(-1,\dots,-1,n,-1,\dots,-1)
\end{equation}
where the different coordinate is the $m$+1-th.
Then
\begin{equation}\label{eq:inp}
<\mu_m,\mu_m> = \frac{n}{n+1} \qquad <\mu_m,\mu_l> = \frac{-1}{n+1} \qquad  <\mu_m+\mu_l, \mu_l-\mu_m> = 0
\end{equation}
for all $m,l$. We didn't need the following two, but we collect them too ($m\neq l$)
$$<\mu_m-\mu_l, \mu_m-\mu_l> = 2 \qquad <\mu_m+\mu_l, \mu_m+\mu_l> = 2\frac{n-1}{n+1}$$
Finally we collect the expression for $\rho$
$$\rho = \sum\omega_i = \frac{n}{2}\,(1,1,\dots,1)-(0,1,\dots,n)$$


\section{{\large Fusion Matrix for the Standard Representation of $U_q(\lie{sl}_{n+1})$}}\label{ap:fus}

We just mention the formula referring to \cite{eo01,ev97,ev99} for the proof and a more detailed development of the theory.
\begin{equation}\label{eq:fus}
J(\lambda) = \sum_{i,j=0}^n E_{i,i}\otimes E_{j,j} + \sum_{i<j} \frac{q-q^{-1}}{1-q^{2(\lambda_i-\lambda_j + j-i)}} \, E_{j,i}\otimes E_{i,j}
\end{equation}    
In our notation $q^{2(\lambda_i-\lambda_j + j-i)} = q^{2<\lambda+\rho,\mu_i-\mu_j>}$. Actually, knowing that $J$ is triangular with $1$s in the diagonal, proposition \ref{prop:sqkz} provides a proof of \eqref{eq:fus}. 


\section{{\large Warnings}}

We list some computational errors we have found in some previous publications. They surely do not affect the theory developed in those papers, but are important if one needs to perform precise calculations.

\begin{enumerate}

\item The version of proposition \ref{prop:factor} found in \cite{fr01} (theorem 4.2 there) is not correct as stated. Also the formula for $r_{V,V}(z)$ on page 28 (which should correspond to $\varrho_{n+1}(z)^{-1}$) is not correct (compare with \eqref{eq:fro}). We used the corrected version found in \cite{efk01} where one can also find a more general expression (for any finite dimensional representation) in the case of $\lie{sl}_2$.

\item The expressions for $\Psi_{\lambda,k}^{v_i\otimes v_j}(x,y)$ and for $R_k(u,\lambda)$ found here can be used to deduce the fusion and exchange matrix without spectral parameter if one picks the corresponding constant coefficients. One then should correct the expression of the exchange matrix of $U_q(\lie{sl}_2)$ for the standard representation found in \cite{eo01} (the term $q^{\frac{n}{n+1}}$ does not appear correctly there).

\item The version of proposition \ref{prop:sqkz} found in \cite{efk01} (chapter 11) does not give fusion solutions. The solutions listed there are actually our solutions \eqref{eq:as}. One should then correct the expression for the connection matrix found in chapter 12. In the case of the standard representation the final answer is given by theorem \ref{thm:formula}. We do not give the correction for the general situation found there since we have not treated it here. We remark that one should read chapter 12 of \cite{efk01} as the guideline for those calculations (as we did here), but should always check the expressions.

\end{enumerate}



\begin{thebibliography}{10}

\bibitem[EFK98]{efk01}
{Etingof, P.I.};~{Frenkel, I.B.}; {Kirillov Jr., A.A.}
\newblock {\em Lectures on representation theory and {K}nizhnik-{Z}amolodchikov
  equations}, volume~58 of {\em Mathematical Surveys and Monographs}.
\newblock AMS, 1998.

\bibitem[ES99]{eo01}
{Etingof, P.}; {Schiffmann, O.}
\newblock Lectures on the dynamical {Y}ang-{B}axter equations.
\newblock {\small QA/9908064}.

\bibitem[EV98]{ev97}
{Etingof, P.}; {Varchenko, A.}
\newblock Solutions of the quantum dynamical {Y}ang-{B}axter equation and
  dynamical quantum groups.
\newblock {\em Comm. Math. Phys}, 196:591--640, 1998.

\bibitem[EV99]{ev99}
{Etingof, P.}; {Varchenko, A.}
\newblock Exchange dynamical quantum groups.
\newblock {\em Comm. Math. Phys}, 205:19--52, 1999.

\bibitem[Fel94]{fel94}
{Felder, G.}
\newblock Elliptic quantum groups.
\newblock {\em Proc. of the ICMP}, Paris, 1994.
\newblock {\small hep-th 9412207}.

\bibitem[FelV99]{felv99}
{Felder, G.}; {Varchenko, A.}
\newblock The elliptic gamma function and
  ${SL}(3,\mathbb{Z})\ltimes\mathbb{Z}^3$.
\newblock {\small QA/9907061}.

\bibitem[FR92]{fr01}
{Frenkel, I.}; {Reshetikhin, N.}
\newblock Quantum affine algebras and holonomic difference equations.
\newblock {\em Comm. Math. Phys.}, 146:1--60, 1992.

\bibitem[GaRa90]{gara01}
{Gasper, G.}; {Rahman, M.}
\newblock {\em Basic hypergeometric series}, volume~35 of {\em Encyclopedia of
  Mathematics and it applications}.
\newblock Cambridge University Press, 1990.

\bibitem[J86]{jim02}
{Jimbo, M.A.}
\newblock A q-analogue of ${U}(\mathfrak{gl}(n+1))$, {H}ecke algebra and the
  {Y}ang-{B}axter equation.
\newblock {\em Lett. Math. Phys.}, 11:257--252, 1986.

\bibitem[TV97]{tv97}
{Tarasov, V.}; {Varchenko, A.}
\newblock Geometry of q-hypergeometric functions, quantum affine algebras and elliptic quantum groups.
\newblock {\em Ast\'erisque}, 246, 1997.
\newblock {\small q-alg/9703044}. 

\bibitem[WW02]{ww01}
{Whittaker, E.T.};{Watson, G.N.}
\newblock {\em A course of modern analysis}.
\newblock Cambridge University Press, 1902.

\end{thebibliography}

\end{document}